\documentclass[a4paper,12pt]{article}

\usepackage{graphicx}      
\usepackage{epsfig} 
\usepackage{amsmath} 
\usepackage{amssymb}  

\gdef \RR{{\mathbb R}}

\gdef \ZZ{{\mathbb Z}}

\newtheorem{theorem}{\bf{Theorem}}

\newtheorem{proposition}[theorem]{\bf{Proposition}}

\title{\bf Positive trigonometric polynomials for strong stability of
difference equations$^1$}

\begin{document}

\footnotetext[1]{This research has been supported by the Ministry of
Education of the Czech Republic under Project 1M0567, by the Grant Agency
of the Czech Republic under Project 103/10/0628, and by a bilateral
Czech-French research project Barrande (MEB020915).}

\author{Didier Henrion$^2$, Tom\'{a}\v{s} Vyhl\'{\i}dal$^3$}

\footnotetext[2]{CNRS; LAAS; 7 avenue du colonel Roche, F-31077 Toulouse, France;
Universit\'e de Toulouse; UPS, INSA, INP, ISAE; LAAS; F-31077 Toulouse, France.
Faculty of Electrical Engineering, Czech Technical University in Prague,
Technick\'a 4, CZ-16626 Prague, Czech Republic, {\tt henrion@laas.fr}}
\footnotetext[3]{Center for Applied Cybernetics, Department of Instrumentation
and Control Engineering, Faculty of Mechanical Engineering, Czech Technical
University in Prague, Technick\'{a} 4, 166 07 Praha 6, Czech Republic,
{\tt tomas.vyhlidal@fs.cvut.cz}}

\maketitle

\begin{abstract} 
We follow a polynomial approach to analyse strong stability of
linear
difference equations with rationally independent delays. Upon application
of the
Hermite stability criterion on the discrete-time homogeneous
characteristic polynomial,
assessing strong stability amounts to deciding positive definiteness of a
multivariate
trigonometric polynomial matrix. This latter problem is addressed with a
converging
hierarchy of linear matrix inequalities (LMIs). Numerical experiments
indicate that
certificates of strong stability can be obtained at a reasonable
computational cost
for state dimension and number of delays not exceeding 4 or 5.
\end{abstract}

\begin{center}
\end{center}


\section{Introduction}

In general, spectrum-based analysis of time-delay systems can be handled in the same way it is done for delay-free systems. Although the spectrum is infinite, stability is determined by the rightmost eigenvalues, more precisely by the sign of the spectral abscissa, the maximum
real part of the eigenvalues. For retarded systems, the spectral abscissa is nonsmooth but continuous in all parameters of the system, including time delays, see \cite{Joris07}. However, it results from \cite{henry, Avellar, Hale93, verduyn}, that, in general, it is not the case for neutral systems and kernel operators - the so-called associated difference equation, see also \cite{wimneutral, contppN, Michiels07}. It is well-known that the spectral abscissa of the difference equation is not continuous in delays. Thus, arbitrarily small changes in the delay values can destroy stability. Moreover, it can even happen that the number of unstable roots increases stepwise from zero to infinity. In order to handle this hypersensitivity of the stability of the difference equation with respect to delay values, the concept of {\em strong stability} was introduced by \cite{verduyn}. Let us remark that the {\em strong stability} concept has recently been generalized by \cite{Wim09} towar
 d difference equations with dependencies in the delays. 

As stability of its kernel operator is a necessary condition for stability of a neutral system, all the hypersensitivity stability issues are carried over to the stability of neutral systems. Thus the {\em strong stability} test should always be performed to guarantee practical stability of neutral systems. However, as will be shown later in the text, the {\em strong stability} test is rather complex. So far, a coarse numerical implementation of the test without guarantee or certificate has been used as a rule, see e.g. \cite{contppN, Tomas10}. Even though this {\em brute force} based approach works in most cases, it might fail due to approximation errors in the numerical scheme. As the main result of this paper we propose a more rigorous {\em strong stability test} that is based on a polynomial approach, relying on the numerical solution of a
hierarchy of linear matrix inequalities (LMIs). 

In the field of time-delay systems, LMIs are usually used as stability determining criteria resulting from the Lyapunov time-domain approach, see e.g. \cite{Silviu} or \cite{Xu-Guang}, among many others.
     
\subsection{Problem statement}
We consider a neutral system of the following form

\begin{equation}\label{sys}
\frac{d}{dt}\left(x(t)+\sum_{k=1}^m H_k x(t-\tau_k)\right)=A_0\ x(t)+\sum_{j=1}^p A_j x(t-\vartheta_j)
\end{equation}

where $x\in\RR^n$ is the state, $\tau_k>0, k=1,\ldots m$ and $\vartheta_j>0, j=1, \ldots p$  are the time delays. It is well-known, see \cite{Hale93}, that a necessary condition for stability of neutral system (\ref{sys}) is stability of the associated difference equation  

\begin{equation}\label{difeq}
x(t)+\sum_{k=1}^m H_k x(t-\tau_k)=0.
\end{equation}

Moreover, {\em strong stability} of equation (\ref{difeq}) is required, i.e. stability independent of the values of the delays,
\cite{Avellar, Hale93}. In \cite{verduyn} (Theorem 2.2 and Corollary 2.2), a condition for strong stability condition is stated
as follows:

\begin{proposition}\label{propstrong} 
Delay difference equation (\ref{difeq}) is strongly stable if and only if

\begin{equation}\label{gamma0}
\gamma_0:=\max_{\theta\in [0,\ 2\pi]^{m}}\
r_{\sigma}\left(\sum_{k=1}^m H_k
e^{-i\theta_k}\right)<1,
\end{equation}

where $r_{\sigma}$ denotes the spectral radius, i.e. the maximum modulus of the eigenvalues.
Furthermore, if $\gamma_0>1$ then equation (\ref{difeq}) is exponentially unstable for rationally independent\footnote{The $m$ numbers $\tau=(\tau_1,\ldots,\tau_m)$ are rationally independent if and only if
$\sum_{k=1}^m n_k\tau_k=0,\ n_k\in\ZZ$
implies $n_k=0,\ \forall k=1,\ldots,m$. For instance, two delays $\tau_1$ and $\tau_2$ are rationally independent if their ratio is an irrational number.} delays.
\end{proposition}

Notice that the quantity $\gamma_0$ does not depend on the value of the delays, i.e. exponential stability locally in the delays is equivalent with exponential stability globally in the delays \cite{verduyn}.

Let us remark that by homogeneity, the expression of $\gamma_0$ can be simplified to
\begin{equation}\label{gamma0simp}
\gamma_0=\max_{\theta\in [0,\ 2\pi]^{m-1}}\
r_{\sigma}\left(\sum_{k=1}^{m-1} H_k
e^{-i\theta_k}+H_m\right).
\end{equation}

We conclude the section with some properties of the quantity $\gamma_0$, see \cite{Michiels07, contppN}, for more details.   

\textbf{Properties}

\begin{enumerate}
	\item Stability of difference equation (\ref{difeq}) with rationally independent delays implies strong stability, and vice versa
	\item In the case of one delay ($m=1$),
	\[
	\gamma_0=r_{\sigma}(H_1).
	\]
\item In the case of a scalar equation ($n=1$),
\[
\gamma_0=\sum_{k=1}^m |H_k|.
\]
\item A sufficient, but as a rule conservative, condition for strong stability is given by
\[
\sum_{k=1}^m \left\|H_k\right\|<1
\]
where $\|.\|$ denotes the matrix Euclidean norm, i.e. the maximum singular value.
\end{enumerate}

\subsection{Computational issues}\label{compissues}
The problem of solving (\ref{gamma0}) can be formulated as an optimization task with the objective to find the global maximum of spectral radius over $\theta\in [0,\ 2\pi]^{m}$. However, in general the objective function $r_{\sigma}(\theta)$ is nonconvex, i.e. it can have multiple local maxima. Besides, the function can be nonsmooth (e.g. at the points where the spectral radius is determined by more than either one single eigenvalue or a couple of complex conjugate eigenvalues). The fact that the function is nonsmooth precludes the use of standard
optimization procedures. Instead, nonsmooth optimization methods can be used, such as gradient sampling, see \cite{OVERTON,Overton09}.
However, even though these methods can handle the problem of nonsmoothness, they converge to local extrema as a rule. As suboptimal solutions
are not sufficient (the global maximum of the spectral radius is needed) a brute force method has been used to solve the task so far,
see \cite{contppN, Wim09, Tomas10}. In the first step, each dimension of $[0,\ 2\pi]^{m}$ is discretized to $N$ points. Then evaluation of (\ref{gamma0}) consists in solving $N^m$ times $n\times n$ eigenvalue problems. Hence, the overall cost of one evaluation of $\gamma_0$ is $O\left(N^{m} n^3\right)$, see \cite{Tomas10}. If the simplified expression (\ref{gamma0simp}), the computational costs reduces to $O\left(N^{m-1} n^3\right)$. Obviously, the complexity of the computation grows considerably with the number of delays in the difference equation. Moreover, the risk of missing global extrema due to sparse or inappropriate gridding cannot be avoided.

\section{Strong stability and Hermite's condition}

Consider the characteristic polynomial 
\begin{equation}\label{polynomial}
p(z) = \det (z_0 I_n + \sum_{k=1}^m z_k H_k),	
\end{equation}
which is homogeneous of degree $n$ in $m+1$
variables $z_k$, $k=0,1,\ldots,m$. 

Based on (\ref{gamma0}), considering $z_k=e^{j\theta_k}, \theta_k\in[0,2\pi], k=1,..,m$, the difference equation (\ref{difeq}) is strongly stable if and only if
the univariate polynomial
\[
z_0 \rightarrow p(z)
\]
is discrete-time stable,i.e. it has all its roots in the open unit disk.

In order to deal with stability of this polynomial, we use
a stability criterion based on the Hermite matrix.
It is a Hermitian matrix of dimension $n$ whose entries are quadratic
in the coefficients of the polynomial. The Hermite
matrix $z_1,\ldots,z_m \rightarrow H(z)$ is therefore a trigonometric
polynomial matrix in $m$ variables $z_1,\ldots,z_m$.

Derived by the French mathematician Charles Hermite in 1854,
the Hermite matrix criterion is a symmetric version of
the Routh-Hurwitz criterion
for assessing stability of
a polynomial. It says that a polynomial $p(z)=p_0+p_1z+\cdots+p_nz^n$
has all its roots in the open upper half of the complex plane if and only if its Hermite
matrix $H(p)$ is positive definite.
Note that $H(p)$ is $n$-by-$n$, Hermitian and
quadratic in coefficients $p_k$, so that the above
necessary and sufficient stability condition is a
quadratic matrix inequality (QMI) in coefficient vector
$p = [p_0\:p_1\:\cdots\: p_n]$.

The standard construction of the Hermite matrix goes through
the notion of B\'ezoutian, a particular form of the resultant.
A bivariate polynomial is constructed, from which a quadratic
term is factored out, yielding a quadratic form shaped by the
Hermite matrix. The construction is explained e.g. in
\cite{hpas03} and references therein. See especially \cite{lbk91}
which explains that a discrete-time Hermite matrix, sometimes called Schur-Cohn or
Schur-Fujiwara matrix, can be obtained similarly.
The discrete-time Hermite matrix is also quadratic in the $p_k$,
and it is positive definite if and only if polynomial $p(z)$ has all
its roots in the open unit disk.

Zden\v ek Hur\'ak pointed out that there is a much simpler construction
of the Hermite matrix in the discrete-time case. The construction
can be traced back to Issai Schur \cite{ps05}, and it is explained in
\cite{hurak}. Entrywise formulas are also described in \cite[Theorem 3.13]{barnett}.
Let
\[
S_1(p) = \left[\begin{array}{cccc}
p_n & p_{n-1} & p_{n-2} \\
0 & p_n & p_{n-1} \\
0 & 0 & p_n \\
& & & \ddots
\end{array}\right] \quad
S_2(p) = \left[\begin{array}{cccc}
p_0 & p_1 & p_2 \\
0 & p_0 & p_2 \\
0 & 0 & p_0 \\
& & & \ddots
\end{array}\right]
\]
be $n$-by-$n$ upper-right triangular Toeplitz matrices. Then
\[
H(p)=S^T_1(p)S_1(p)-S^T_2(p)S_2(p).
\]

Strong stability of the difference equation is hence equivalent
to positive definiteness of the Hermite matrix of the univariate
characteristic polynomial, which is a multivariate trigonometric
polynomial matrix in $z_1,\ldots,z_m$. We express this constraint as
\begin{equation}\label{hermpos}
H(z_1,\ldots,z_m) \succ 0.
\end{equation}

\section{Positivity of trigonometric polynomials}

As shown in the previous section, the key ingredient in our approach
to strong stability of difference equation is assessing
positivity of multivariate trigonometric polynomials.
This topic has been subject to recent studies, and the
recent monograph \cite{dumitrescu} is a good introduction
focusing on signal processing applications. 

In this section we start with a scalar multivariate trigonometric polynomial,
formulate its positivity test as a minimization problem,
describe an LMI hierarchy yielding an asymptotically
converging monotonically increasing sequence of lower bounds.
We also describe a hierarchy of eigenvalue problems (linear algebra,
much simpler computationally that LMI methods)
to generate a hierarchy of upper bounds.

Then we extend these results to matrix polynomials, and
describe the hierarchy of LMI problems that must be solved
to guarantee positivity of a trigonometric matrix polynomial
at the price of solving a hierarchy of convex problems,
the decision variables being entries of a Gram matrix
yielding a sum-of-squares decomposition for the
matrix polynomial.

\subsection{Minimising trigonometric polynomials}

A trigonometric polynomial has the form
$h(z) = \sum_{\alpha} h_{\alpha} z^{\alpha}$
where integer vector $\alpha \in {\mathbb N}^n$ is a multi-index such that
$z^{\alpha} = \prod_{i=1}^n z_i^{\alpha_i}$, complex vector
$z \in {\mathbb C}^n$ contains indeterminates such that $z_i = e^{j\theta_i}$
for some $\theta \in [0,2\pi]^n$, and complex numbers
$h_{\alpha} \in {\mathbb C}$ are coefficients.
We use the notation $z \in {\mathbb T}^n$ to
capture the constraint that each variable $z_k \in {\mathbb C}$
belongs to the unit disk ${\mathbb T}$.

We consider real trigonometric polynomials such that $h(z) = h(z)^*$
where the star denotes complex conjugation. These are such that
$\sum_{\alpha} h_{\alpha} z^{\alpha} = \sum_{\alpha} h_{\alpha}^* z^{-\alpha}$
and hence $h_{\alpha} = h_{-\alpha}^*$.

Since $h(z)$ maps ${\mathbb T}^n$ onto $\mathbb R$,
we are interested in solving the problem
\[
h_{\min} = \min_{z \in {\mathbb T}^n} h(z).
\]

\subsection{Hierarchy of lower bounds via SDP}\label{sdp}

In this section we construct a monotonically decreasing
sequence of lower bounds on $h_{\min}$ that converges asymptotically.
Each bound can be computed by solving an LMI, a
convex semidefinite programming (SDP) problem.

First note that by definition
\begin{equation}\label{meas}
h_{\min} = \min_{\mu} \int_{{\mathbb T}^n} h(z) d\mu(z)
\end{equation}
where the minimisation is over all probability measures
defined on the sigma-algebra of the multidisk ${\mathbb T}^n$, see Chapter 5 in \cite{lasserre}.

Let us express polynomial $h(z)$ as a Hermitian quadratic form
\begin{equation}\label{sos}
h(z) = b_k^*(z){\mathbf X}_k b_k(z)
\end{equation}
where $b_k(z)$ is a vector basis of trigonometric polynomials
of degree up to $k$, e.g. containing monomials $z^{\alpha}$,
$\alpha \geq 0$, $\max_{i=1,\ldots,m} \alpha_i \leq k$.
Matrix ${\mathbf X}_k$ is called the
Gram matrix of polynomial $h(z)$ in basis $b_k(z)$.
Then a result of functional analysis by M. Putinar, transposed
to trigonometric polynomials \cite[Theorems 3.5 and 4.11]{dumitrescu},
states that
$h(z) > 0$ if and only if there exists a finite integer $d$
and a positive semidefinite Hermitian matrix $X_d \succeq 0$
such that (\ref{sos}) holds for $k=d$. 

As soon as $k$ is fixed, finding a matrix ${\mathbf X}_k \succeq 0$
satisfying (\ref{sos}) can be cast into an SDP 
feasibility problem which amounts to
expressing polynomial $h(z)$ as
a sum-of-squares (SOS) of trigonometric polynomials
of degree $k$.

Now defining
\[
\begin{array}{rcl}
\underline{h}_k = & \sup & \underline{h} \\
& \mathrm{s.t.} & h(z)-\underline{h} = b_k^*(z) {\mathbf X}_k b_k(z)
\;\mathrm{for}\;\mathrm{some}\: {\mathbf X}_k \succeq 0
\end{array}
\]
it follows that $\underline{h}_k \leq \underline{h}_{k+1}$
and we expect that $\lim_{k \rightarrow \infty} \underline{h}_k = h_{\min}$,
even though a rigorous proof of convergence is out of the
scope of this paper.

\subsection{Hierarchy of upper bounds via EVP}\label{evp}

In this section we show that we can construct a monotonically increasing
sequence of upper bounds on $h_{\min}$ that converges asymptotically.
Each bound can be computed by solving an eigenvalue problem (EVP)

In problem (\ref{meas}) let us consider that measure $\mu$
is absolutely continuous w.r.t. measure $\nu$,
the probability measure supported uniformly on the multidisk.
Let us further restrict the class of measures by considering
that there exists a trigonometric polynomial $q_k(z) = \sum_{0 \leq \alpha \leq k} {q_k}_{\alpha} z^{\alpha} = {\mathbf q}_k^* b_k(z)$
of total degree $k$
such that $\mu_k(dz) = q_k^*(dz)q_k(dz)\nu(dz)$, with
$\lim_{k \rightarrow \infty} \mu_k = \mu$ since ${\mathbb T}^n$
is compact. Let $y_{\alpha} = \int_{{\mathbb T}^n} z^{\alpha}d\nu(z)$ denote the moment of order $\alpha$ of $\nu$. Finally, let us define
\[
\overline{h}_k = \min_{\mu_k} \int_{{\mathbb T}^n} h(z) d\mu_k(z)
\]
as an optimisation problem over this restricted class of measures.

With these notations
\[
\int h(z)d\mu_k(z) = \int h(z)q_k^*(z)q_k(z)d\nu(z) =
\]
\[
=\int h(z){\mathbf q}_k^* b_k(z)b_k^*(z){\mathbf q}_k d\nu(z)
\]
is the same as
\[
{\mathbf q}_k^* \left(\int h(z)b_k(z)b_k^*(z)d\nu(z) \right) {\mathbf q}_k =
{\mathbf q}_k^* {\mathbf M}_k(h\:y) {\mathbf q}_k
\]
where ${\mathbf M}_k(h\:y)$ is called the localising matrix of order $k$
of measure $\nu$ w.r.t. polynomial $h$, see \cite{lasserre}.
Its rows and columns are indexed by multi-indices $\beta$ and $\gamma$
respectively, and its entry $(\beta,\gamma)$ is equal to
$\sum_{\alpha} h_{\alpha} y_{\alpha-\beta+\gamma}$. Therefore
matrix ${\mathbf M}_k(h\:y)$ can be obtained from the moments of $\nu$,
and hence it is given. It is positive definite.

If $h(z)=1$, matrix ${\mathbf M}_k(y)$
is called the moment matrix of order $k$ of
measure $\nu$. Its entry $(\beta,\gamma)$ is equal to $y_{-\beta+\gamma}$,
and hence matrix ${\mathbf M}_k(y)$ is given as well. Since $\mu_k$ is
a probability measure
\[
\int d\mu_k = \int q_k^*q_k d\nu_k =
{\mathbf q}_k^* {\mathbf M}_k(y) {\mathbf q}_k = 1
\]
and hence 
\[
\begin{array}{rcl}
\overline{h}_k = & \min_{{\mathbf q}_k} & {\mathbf q}_k^*
{\mathbf M}_k(h\:y)
{\mathbf q}_k \\
& \mathrm{s.t.} & {\mathbf q}_k^* {\mathbf M}_k(y) {\mathbf q}_k = 1.
\end{array}
\]
It follows that $\overline{h}_k \leq \overline{h}_{k+1}$
and $\lim_{k \mapsto \infty} \overline{h}_k = h_{\min}$
even though I am not totally confident that 
this latter result is correct.

Finally,
given positive definite Hermitian matrices $A$ and $B$,
optimisation problem $\min_v v^*Av$ s.t. $v^*Bv = 1$
can be solved via linear algebra. Indeed, let $z$
denote an eigenvalue of the pencil $zB-A$, and let $\bar{v}$
denote the corresponding unit eigenvector. Then vector
$v=(\bar{v}^*B\bar{v})^{-\frac{1}{2}}\bar{v}$
is such that $v^*Bv=1$
and $v^*Av=z$. Minimising this quantity
then amounts to finding the minimum eigenvalue of
pencil $zB-A$.

\subsection{Polynomial matrices}\label{polymat}

The above results on scalar polynomials can be extended
directly to polynomial matrices by considering a
matrix basis instead of a vector basis to build
the Hermitian matrix representation (\ref{sos}).

In the context of our strong stability analysis problem,
the core idea is then to replace the (typically difficult)
Hermite matrix positivity condition (\ref{hermpos})
with a hierarchy of tractable SDP problems.
We write
\begin{equation}\label{matsos}
\begin{array}{rcl}
\underline{h}_k = & \sup & \underline{h} \\
& \mathrm{s.t.} & H(z)-\underline{h} = (b_k(z) \otimes I_n)^* {\mathbf X}_k (b_k(z) \otimes I_n)\\
& & \quad {\mathbf X}_k \succeq 0
\end{array}
\end{equation}
as an LMI relaxation of order $k$ of positivity
condition (\ref{hermpos}).

If $\underline{h}_k > 0$
for some $k$, then it implies that (\ref{hermpos})
is satisfied.

If $\underline{h}_k \leq 0$ for some $k$, then
we cannot conclude directly, but we can try to
extract from the dual (moment) SDP problem
a certificate that indeed matrix $H(z)$ cannot
be positive definite, see \cite{tac06} even though
the trigonometric polynomial matrix case is not
developed in this reference. If we cannot extract
useful information from the dual problem, we have
to increase the value of $k$ and solve the
next LMI in the hierarchy.

\section{Complexity}

Let $M$ denote the size of the Gram matrix ${\mathbf X}_k$ in SDP problem (\ref{matsos}).
If we use an interior-point method, the worst-case complexity of one Newton
iteration for an SDP problem in a cone of that size is $O(M^6)$. Experiments
reveals that the practical complexity is approximately $O(M^4)$.

The number of monomials of $m$ variables of degree $k$ in basis $b_k(z)$
is equal to $(k+1)^m$. Polynomial $p(z)$ has $m$ variables and degree $n$
so degree $k$ in (\ref{sos}) should be such that $2k \geq n$. Note that we
can have $2k > n$ since higher-degree terms may cancel in the right handside
of equation (\ref{sos}).

If we choose $k = n/2$ or $k = (n+1)/2$ depending on whether $n$ is even or not,
in terms of complexity $M=O(n^{m+1})$. The overall complexity of our SDP
approach to strong stability analysis therefore grows exponentially in the number
of delays $m$, and polynomially in the number of states $n$. However the exponent
of this polynomial growth is quite large. In comparison, the gridding approach
mentioned at the beginning of the paper has a complexity which also grows
exponentially in the number of delays, but the dependence on the number
of states is only cubic. However, contrary to the SDP approach,
the gridding approach does not provide guarantees.

\section{Examples}

Preliminary numerical examples indicate that the
EVP approach of paragraph \ref{evp} yields
a sequence of bounds which converges slowly (sublinearly).
This is why in this section we focus exclusively on
the SDP approach of paragraph \ref{sdp}.

We implemented a collection of Matlab functions
to manipulate trigonometric polynomials,
Hermite matrices, and formulate SDP problems
corresponding to positivity checks. The functions
are available for download\footnote{
\tt homepages.laas.fr/henrion/software/trigopoly.tar.gz}
and they provide the following functionalities:
\begin{itemize}
\item {\tt sampledet.m} - given a collection of
matrices $H_k$, $k=0,\ldots,m$, this function
computes the coefficients of the multivariate
polynomial $p(z)=\det(H_0+H_1z_1+\cdots+H_mz_m)$;
it proceeds by sampling and interpolation,
as described in \cite{cdc05}
\item {\tt trigoherm.m} - computes the Hermite matrix
of a homogenized multivariate polynomial; it uses
the formula of \cite[Theorem 3.13]{barnett} adapted to
complex coefficients
\item {\tt trigohermgram.m} - computes the SDP problem
corresponding to the positivity test for a given
Hermitian multivariate polynomial matrix; the SDP problem
is given in SeDuMi's input format
\[
\begin{array}{llll}
\min & c^T x & \ \ \ \max & b^T y \\
\mathrm{s.t.} & Ax = b & \ \ \ \mathrm{s.t.} & z = c-A^Ty \\
& x \in K & & z \in K
\end{array}
\]
where $x \in {\mathbb R}^N$, $y \in {\mathbb R}^M$,
and $K$ is the cone of positive semidefinite matrices
of size $S=\sqrt{N}$.
\end{itemize}
Some instrumental functions are also provided,
namely {\tt genmon.m} which generates powers of monomials
and {\tt locmon.m} which locates a monomial in a Gram
matrix. Besides, the function {\tt bfssde.m} is available to evaluate (\ref{gamma0}) by brute force, as explained in subsection \ref{compissues}

\subsection{Three states, two delays}\label{n3m2}
\begin{figure}
	\centering
		\includegraphics[width=0.7\textwidth]{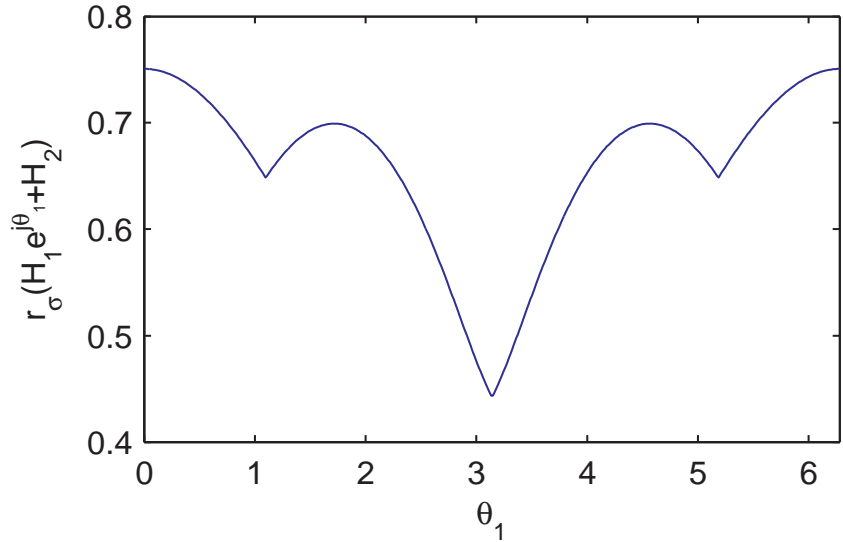}
	\caption{Spectral radius $r_\sigma(\theta_1)$ for the example in Subsection \ref{n3m2}}
	\label{spectrradius1}
\end{figure}

We adopt the illustrative example from \cite{contppN} with $n=3$, $m=2$, where  
\[
H_1 = \left[\begin{array}{rrr}
0 & 0.2 & -0.4 \\
-0.5 & 0.3 & 0 \\
0.2 & 0.7 & 0
\end{array}\right], \quad
H_2 = \left[\begin{array}{rrr}
-0.3 & -0.1 & 0 \\
0 & 0.2 & 0 \\
0.1 & 0 & 0.4
\end{array}\right]
\]
for which {\tt bfssde.m} (with $N=360$) provides $\gamma_0=0.7507$ in less then $0.1$ seconds
under Matlab 7.7 on our Linux PC
equipped with Intel Xeon 2.67GHz CPU with 8GB RAM. On 
Fig. \ref{spectrradius1} shows the spectral radius as a function of $\theta_1$.

The following Matlab script assesses stability
of the corresponding difference equation
by first building the determinantal polynomial,
then the corresponding Hermite matrix, then
the SDP problem, and eventually by solving
the SDP problem with SeDuMi,
a primal-dual interior-point solver:
\small
\begin{verbatim}
H1=[0 0.2 -0.4;-0.5 0.3 0;0.2 0.7 0];
H2=[-0.3 -0.1 0;0 0.2 0;0.1 0 0.4];
p=sampledet({eye(3),H1,H2}); % evaluate determinant
p=p(:,abs(p(1,:))>1e-8); % remove small coefficients
H=trigoherm(p); % compute Hermite matrix
[A,b,c,K]=trigohermgram(H); % build SDP problem
[x,y,info]=sedumi(A,b,c,K); % solve SDP problem
\end{verbatim} 
\normalsize
The resulting SDP problem has size $N=2304$,
$M=225$ and a positive semidefinite Gram matrix 
of size $S=48$ is found after less than 0.1 seconds
with SeDuMi 1.3.

We can also specify the strong stability radius
$\gamma_0$ as a second input argument to function {\tt trigoherm}.
Internally, the polynomial is scaled appropriately
and positivity of the Hermite matrix is assessed:

\small
\begin{verbatim}
H=trigoherm(p,0.750);
[A,b,c,K]=trigohermgram(H); 
[x,y,info]=sedumi(A,b,c,K); 
\end{verbatim}
\normalsize

With the above sequence the SDP problem is found
feasible. Changing the first instruction to

\small
\begin{verbatim}
H=trigoherm(p,0.751);
\end{verbatim}
\normalsize

makes the resulting SDP problem infeasible,
and this is certified by SeDuMi which returns
a dual Farkas vector. As discussed at the end
of paragraph \ref{polymat}, further analysis
is required to conclude that indeed the Hermite
matrix cannot be positive definite. We leave
a comprehensive treatment of this case
for further work.

\subsection{Four states, three delays}\label{n4m3}
We consider a system with $n=4$, $m=3$, where

\begin{figure}
	\centering
		\includegraphics[width=0.7\textwidth]{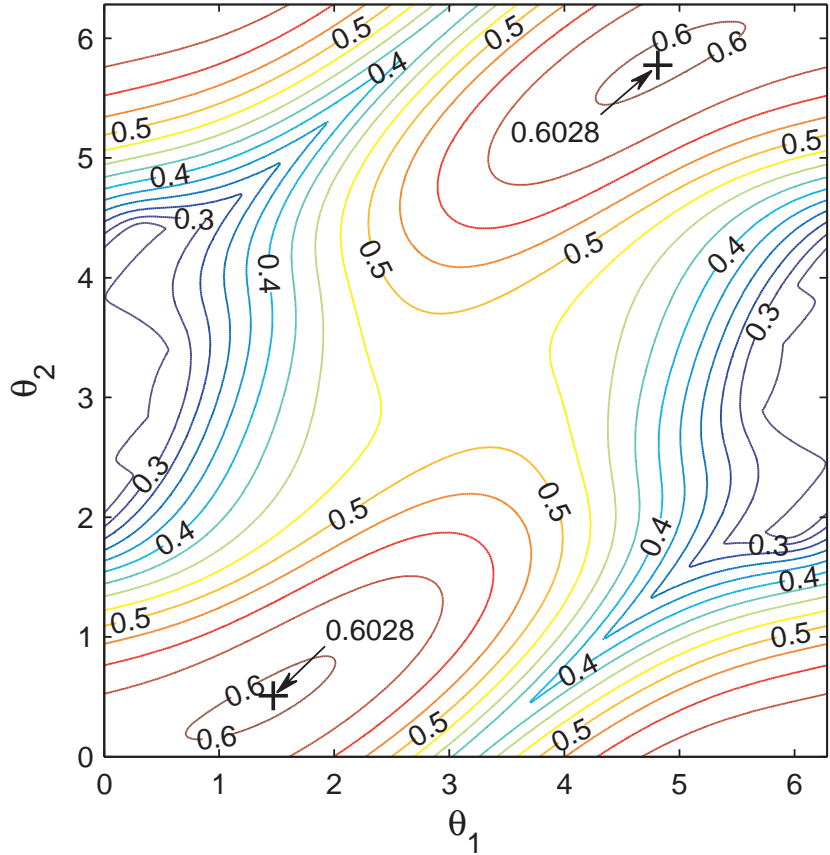}
	\caption{Spectral radius $r_\sigma(\theta_1,\theta_2)$ for the example in Subsection \ref{n4m3}}
	\label{spectrradius2}
\end{figure}

\small
\begin{verbatim}
H1=[-0.15 0 0.32 0;0 -0.07 0 0.05;
    0.08 0 0.04 0;0.2 0.03 0 -0.13];
H2=[-0.02 0.12 0 0.25;0 -0.05 0.04 0;
    0 0.23 0 -0.3;0.19 0 0.28 -0.09];
H3=[0 0 -0.03 0.14;0.01 -0.04 0 0;
    0 0 0.09 0.26; 0.05 -0.27 -0.06 0];
\end{verbatim}
\normalsize

for which {\tt bfssde.m} (with $N=360$) provides $\gamma_0=0.6028$ in 4.5 seconds,
see Fig. \ref{spectrradius2} with the distribution of the spectral radius with respect to values of $\theta_1$ and $\theta_2$.

The resulting SDP problem has size $N=250000$,
$M=5840$ and a positive semidefinite Gram matrix 
of size $S=500$ is found after approximately 6 minutes
of CPU time, certifying that the spectral radius is less
than one.

\subsection{Four states, four delays}

We conclude with an example with $n=4$, $m=4$ and the matrices

\small
\begin{verbatim}
H1=[0.1 0 0 -0.2;pi/5 -0.1 0 -0.3;
    0 0 0.03 2;0 -exp(-1) 0 0.23]
H2=[0 0 0 0.0456;0 -0.33 0.11 0;
    0 1 0.2 0;0 -exp(-3) 0.176 0.73]
H3=[0.1 0.65 0 0.42;0.087 -pi/8 -0.1 0;
    0 -0.063 0 0.72;0.076 0.1 0 -0.23]
H4=[-0.678 0 0 -0.4;-0.0983 0 0 0;
    0 0.0763 0 0.2;-exp(-5) 0 0.36 0]
\end{verbatim}
\normalsize

for which {\tt bfssde.m} (with $N=360$) provides $\gamma_0=1.7649$ in more than
30 minutes.

%
The resulting SDP problem has size $N=6250000$,
$M=52496$ and a positive semidefinite Gram matrix 
of size $S=2500$. This problem cannot not be solved
on our computer, SeDuMi issues an out of memory error message.
In this case, we may to try to exploit the problem structure
(sparsity) to generate a smaller SDP problem, but this is
out of the scope of this paper.

\section{Conclusions}

In the context of neutral time-delay systems,
strong stability of difference equations is generally
assessed numerically with a brute force gridding approach.
A parallel can be draw with the $\mu$-analysis approach
to robustness of linear systems, see e.g. \cite{doyle}
where brute force gridding can yield misleading results
and should be replaced, if possible, with more rigorous
certificates of robustness.

In this paper, using the Hermite stability criterion for discrete-time
polynomials the problem of assessing strong stability is reformulated
as the problem of deciding positive definiteness of a
trigonometric matrix polynomial of size equal to the state dimension
and number of variables equal to the number of delays. This
decision problem is hard, but it can be approached through a
converging hierarchy of tractable semidefinite programming (SDP) or
linear matrix inequality (LMI) relaxations. Numerical experiments
reveal that the approach is limited to small state dimension
and a small number of delays, as expected.

\section*{Acknowledgements}

This work benefited from exchanges with
Jean-Jacques Loiseau and Wim Michiels.

\end{document}